\documentclass{lmcs}
\pdfoutput=1

\usepackage[utf8]{inputenc}
\usepackage{lastpage}
\lmcsdoi{16}{3}{5}
\lmcsheading{}{\pageref{LastPage}}{}{}%
{Jan.~10,~2017}{Jul.~14,~2020}{}


\ACMCCS{[{\bf Theory of computation}]: Logic---Constructive Mathematics; [{\bf Mathematics of
  computing}]: Continuous Mathematics---Continuous Functions}

\usepackage{hyperref}
\hypersetup{hidelinks}
\usepackage{tikz}

\newenvironment{ProofOf}[1][\misc]{\par
  \trivlist
  \item[\hskip\labelsep
        \itshape
  Proof of #1]\ignorespaces
}
{\endtrivlist}

\def\ie{{\em i.e.}}

\begin{document}

\title{Interpolating Between Choices for the Approximate Intermediate Value Theorem}

\author[Matthew Frank]{Matthew Frank}	
\address{Independent, Astoria, NY}	
\email{mfrank@uchicago.edu} 
\thanks{I thank Andrej Bauer, Mike Shulman and Bas Spitters for helpful dialogue on this topic.}	

\begin{abstract}
  \noindent This paper proves the approximate intermediate value theorem, constructively and from notably weak hypotheses:  from pointwise rather than uniform continuity, without assuming that reals are presented with rational approximants, and without using countable choice.  The theorem is that if a pointwise continuous function has both a negative and a positive value, then it has values arbitrarily close to 0.  The proof builds on the usual classical proof by bisection, which repeatedly selects the left or right half of an interval; the algorithm here selects an interval of half the size in a continuous way, interpolating between those two possibilities.

\end{abstract}

\maketitle

\section{Introduction}\label{S:one}

On MathOverflow [Shu16], Mike Shulman asked:  ``Can the approximate Intermediate Value Theorem be proven constructively about an arbitrary (pointwise) continuous function $f$, without using any form of choice or excluded middle (e.g.\ in the mathematics valid in any elementary topos with a natural numbers object)?"

Here, the approximate Intermediate Value Theorem is
\begin{thm}
If $f$ is pointwise continuous on $[a,b]$ and $f(a)<0<f(b)$, then for every $\epsilon>0$, there is some $x$ with $|f(x)|<\epsilon$. 
\end{thm}

The proof below answers ``Yes" to Shulman's question, satisfying four notable constraints.  The first is that the proof should be constructive, in particular using only intuitionistic logic and not relying on the principle of excluded middle.  This puts the standard intermediate value theorem, in which $f(x)=0$, out of reach [BR87, chapter 6.2].

The second constraint is that the proof should use pointwise rather than uniform continuity.  Under uniform continuity, there is a $\delta$ corresponding to the $\epsilon$, there is a sequence of $x$'s going from $a$ to $b$ in steps of less than $\delta$, and one of those $x$'s will have $|f(x)|<\epsilon$.  This is the proof in Bishop and Bridges [BB85, theorem 4.8, p.\ 40], but pointwise continuity does not constructively imply uniform continuity.

The third constraint is that the proof should not use countable choice.  Under countable choice, it suffices in the proof below to let $d_n=0$ if $f(c_n)<\epsilon$ or $d_n=1$ if $f(c_n)>-\epsilon$, but without countable choice, such a sequence is ill-defined.

The fourth constraint is that the proof should not assume that reals are presented with rational approximants.  If a real number is defined as a sequence of rationals, then $d_n$ can be defined as 0 if the $n^{th}$ rational approximation to $f(c_n)$ is negative, or 1 otherwise.  However, this operation does not respect the relation of equality on reals, and it is not appropriate for reals defined as Dedekind cuts.

The ``yes" is qualified only by the need for real numbers to be defined so that Cauchy sequences of them converge in the topos.  For instance, this result holds if reals are defined in the usual constructive way, as Cauchy sequences of rationals with specified modulus (\ie\ sequences $q_m$ where $|q_m-q_{m'}|<1/m + 1/{m'}$); or if reals are defined as Dedekind cuts [BR87, chapters 2.2 and 7.3].  However, this result may not hold if reals are defined as naive Cauchy sequences of rationals (\ie\ sequences $q_m$ where for all $\epsilon$, for some $n$, for all $m>n$ and $m'>n$, $|q_m-q_{m'}|<\epsilon$), or if they are defined as equivalence classes [Sch03, notes 12-13].

The proof below satisfies all the constraints above, and works for topoi in which countable choice fails.  It strengthens some of Schuster's results [Sch03, section 5].  It can also apply to the case $f(a)\le 0$, $f(b)\ge 0$ by changing strict to non-strict inequalities appropriately.

This proof reduces to the usual classical proof by bisection if the definition of $d_n$ is replaced by $d_n = 0$ if $f(c_n)<0$ and $d_n=1$ otherwise, choosing either the left half or the right half of the interval at each step.  The key difference from that usual proof is that the algorithm here selects an interval of half the size in a continuous way, interpolating between those two possibilities.

\section{Proof of the Theorem}

Define the following inductively:
\[a_1 = a\]
\[b_1 = b\]
\[c_n = (a_n+b_n) / 2\]
\[d_n = \max( 0, \min( \textstyle\frac{1}{2}+ \frac{ f(c_n)}{\epsilon}, 1))\]
\[a_{n+1} = c_n - d_n (b-a)/2^n \]
\[b_{n+1} = b_n - d_n (b-a)/2^n\]
Then $b_n - a_n = (b-a)/2^{n-1}.$

So the $c_n$'s converge to some $c.$  In more detail:  The $d_n$'s are between 0 and 1.  Each interval $[a_{n+1},b_{n+1}]$ is of length $(b-a)/2^n$, and contained in the previous $[a_n, b_n]$.  The $c_n$ form a Cauchy sequence, since for $m<n$, $|c_m-c_n| \le  (b-a)/2^{m-1}$.  

The following diagrams illustrate how an interval $[a_n,b_n]$ and the value $f(c_n)$ determine the next interval $[a_{n+1},b_{n+1}]$.
High values of $f(c_n)$ push the next interval to the left, and low values of $f(c_n)$ push the next interval to the right.

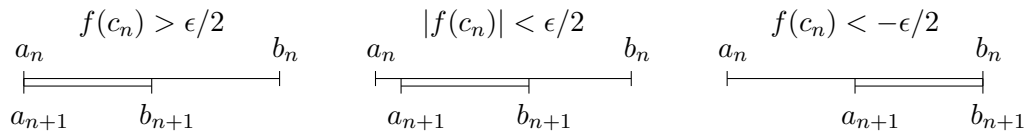
\begin{tikzpicture}[xscale=0.85]

\newcommand{\abc}[4]
{
\node[below] at (#1,-0.2) {$\ \ \phantom{b}a_{n+1}$};
\node[below] at (#2,-0.2) {$\ \ \phantom{b}b_{n+1}$};
\draw (#1,-0.1) -- (#2,-0.1);
\draw (#1,0) -- (#1,-0.2);
\draw (#2,0) -- (#2,-0.2);
\draw (#3-2,0.1) -- (#3-2,-0.1);
\draw (#3+2,0.1) -- (#3+2,-0.1);
\draw (#3-2,0) -- (#3+2,0);
\node[above] at (#3-2,0.1) {$\phantom{b}a_n$};
\node[above] at (#3+2,0.1) {$\phantom{b}b_n$};
}

\node[above] at (1.0,0.4) {$f(c_n)>\epsilon/2$};
\node[above] at (6.5,0.4) {$|f(c_n)|<\epsilon/2$}; 
\node[above] at (12,0.4) {$f(c_n)<-\epsilon/2$};
\abc{-1.0}{1.0}{1.0}{1.6};
\abc{4.9}{6.9}{6.5}{0.75};
\abc{12}{14}{12}{-1.6};

\end{tikzpicture}

\begin{clm}
For any $m\in\mathbb{N}$, either (i) $\exists j \le m,\ |f(c_j)| < \epsilon$
or (ii) $f(a_m) < 0 $ and$ f(b_m) > 0$. 
\end{clm}

\begin{ProofOf}[Theorem from Claim.]
By pointwise continuity at $c$,
let $\delta$ be such that $ |x-c|<\delta$
implies $|f(x)-f(c)| < \epsilon.$
Choose $m$ such that $|c-c_m| < \delta / 2$ and $(b-a)/2^m < \delta / 2$, and apply the claim.

In the first case of the claim, the theorem is immediate.

In the second case of the claim, 
\[|c-a_m| \le |c-c_m| + |c_m-a_m| < \delta, \text{ so }|f(c)-f(a_m)| < \epsilon\]
\[|c-b_m| \le |c-c_m| + |c_m-b_m| < \delta, \text{ so }|f(c)-f(b_m)| < \epsilon\]
So $f(c)$ is within $\epsilon$ of both a negative and a positive number, and $|
f(c)| < \epsilon$.
\end{ProofOf}

\begin{ProofOf}[Claim, by induction on $m$.]
The base case is given by $ f(a) < 0, \ f(b) > 0.$

Now assume the claim for $m$. In case (i), for some $ j < m,\ |f(c_j)| < \epsilon$, and the inductive step is trivial.

In case (ii), use trichotomy with either
$f(c_m) < -\epsilon/2, \ |f(c_m)| < \epsilon$, or $f(c_m) > \epsilon/2.$

If $|f(c_m)| < \epsilon$, then the inductive step is again trivial.

If $f(c_m) > \epsilon/2$, then
\[d_m = 1\]
\[a_{m+1} = a_m, \text{ so }f(a_{m+1}) < 0\]
\[b_{m+1} = c_m,\text{ so } f(b_{m+1}) > 0\]

If $f(c_m) < - \epsilon/2$, then
\[d_m = 0\]
\[a_{m+1} = c_m, \text{ so }f(a_{m+1}) < 0\]
\[b_{m+1} = b_m,\text{ so } f(b_{m+1}) > 0\]

This establishes the claim and the theorem.  \qed 
\end{ProofOf}

\section{Comments and Examples}

Note that $|f(c)|$ may be small or large.  This is especially visible in the following diagrams, where $f$ is piecewise polynomial, the black dots show points $(c_n, f(c_n))$, and the white squares show $(c,f(c))$.

The example below left starts from $[a,b]=[-1,1]$,
$f(x)=\min(\frac{1+6x^2}{7},8+9x)$  and $(c_1,f(c_1))=(0,1/7)$.  Using $\epsilon=1/3$,
the algorithm converges on $c=-0.8894$.... In this case $c$ is not quite the exact zero at $-8/9$, but the difference is below the scale of the graph, and $|f(c)|$ is well under $\epsilon$.  

The example below right starts from the same initial data but with $\epsilon=1/2$.  Then
the algorithm converges on $c=-0.7485$..., for which $f(c) > \epsilon$.  
Taking $\delta=1/8$, $m = 6$, the proof guarantees only that at least one of $f(c_1), \ldots, f(c_6), f(c)$ is within $\epsilon$ of zero, which is enough to show the approximate Intermediate Value Theorem.

\hspace*{-2em}
\begin{tikzpicture}[xscale=3, yscale=2][

\draw[fill] (0,1.25) circle [radius=0.01pt];

\newcommand{\h}[3]
{
\draw (#1, #3) -- (#2, #3);
\draw (#1, #3 + 0.07) -- (#1, #3 - 0.07);
\draw (#2, #3 + 0.07) -- (#2, #3 - 0.07);
}

\newcommand{\axescurvedot}[1]
{
\h{#1-1}{#1+1}{0};
\draw (#1,-1) -- (#1,1);
\node[left] at (#1,1) {$1$};
\node[left] at (#1,-1) {$-1$};
\node[below] at (#1-1,-0.07) {$-1\phantom{-}$};
\node[below] at (#1+1,-0.07) {$1$};
\draw[domain=#1-1:#1-0.810] plot (\x, {8+9*(\x-#1)});
\draw[domain=#1-0.810:#1+1] plot (\x, {(1+6*(\x-#1)*(\x-#1))/7});
\draw[fill] (#1+0,1/7) circle [radius=.50pt];
}

\axescurvedot{-0.1}
\draw[fill] (-0.1+0,0.143) circle [radius=.46pt];
\draw[fill] (-0.1+-0.429,0.300) circle [radius=.46pt];
\draw[fill] (-0.1+-0.679,0.538) circle [radius=.46pt];
\draw[fill] (-0.1+-0.804,0.696) circle [radius=.46pt];
\draw[fill] (-0.1+-0.866,0.205) circle [radius=.46pt];
\draw[fill] (-0.1+-0.897,-0.076) circle [radius=.46pt];
\draw[fill=white] (-0.1+-0.889-0.02,-0.004-0.03) rectangle (-0.1+-0.889+0.02,-0.004+0.03);

\draw[dotted] (-0.1+-1,1/3) -- (-0.1+1,1/3);
\node[left] at (-0.1+-1,1/3) {$\epsilon=\frac{1}{3}$};

\axescurvedot{2.5}
\draw[fill] (2.5+-2/7,73/343) circle [radius=.46pt];
\draw[fill] (2.5+-0.499,0.356) circle [radius=.46pt];
\draw[fill] (2.5+-0.624,0.476) circle [radius=.46pt];
\draw[fill] (2.5+-0.686,0.546) circle [radius=.46pt];
\draw[fill] (2.5+-0.717,0.584) circle [radius=.46pt];
\draw[fill=white] (2.5+-0.749-0.02,0.623-0.03) rectangle (2.5+-0.749+0.02,0.623+0.03) ;

\draw[dotted] (2.5+-1,1/2) -- (2.5+1,1/2);
\node[left] at (2.5+-1,1/2) {$\epsilon=\frac{1}{2}$};
\end{tikzpicture}


\end{document}